\documentclass{article}

\usepackage{amssymb}
\usepackage{amsthm}

\newtheorem{lem}{Lemma}[section]

\newtheorem{theo}[lem]{Theorem}

\newtheorem{fact}[lem]{Fact}
\newtheorem{cor}[lem]{Corollary}

\newtheorem*{theorem}{Main Theorem}
\newtheorem*{definition}{Definition}

\newcommand{\SL}{\mathop{\rm SL}}
\newcommand{\PSL}{\mathop{\rm PSL}}

\newcommand{\<}{\langle}

\renewcommand{\>}{\rangle}
\renewcommand{\o}{^\circ}

\newcommand{\Z}{\mathbb{Z}}
\newcommand{\N}{\mathbb{N}}

\title{Lie rank in groups of finite Morley rank with solvable local subgroups}
\author{Adrien Deloro\footnote{Institut de Math\'ematiques de Jussieu,
UMR 7586, Universit\'e Pierre et Marie Curie (Paris 6), 4 place Jussieu, case 247,
75252, Paris Cedex 05, France} and \'Eric Jaligot\footnote{Institut Fourier, CNRS, Universit\'e Grenoble I, 100 rue des maths, BP 74, 38402 St Martin d'H\`eres cedex, France}}
\date{\today}

\begin{document}
\maketitle

\begin{abstract}
We prove a general dichotomy theorem for groups of finite Morley rank with solvable local subgroups and
of Pr\"{u}fer $p$-rank at least $2$, leading either to some $p$-strong embedding, or to the Pr\"{u}fer $p$-rank being exactly $2$.
\end{abstract}



\vspace{1cm}
\hrule\begin{center}
\'Eric Jaligot left us on July 10, 2013.
\end{center}
\hrule
\vspace{1cm}

\section{Introduction}

The ``size'' of a simple algebraic group over an algebraically closed field can be captured by several means. 
One can measure its Zariski dimension, but one can also consider its Lie rank, 
which is the Zariski dimension of its maximal algebraic tori.
For instance, it is often straightforward to argue by induction on the Zariski dimension, 
as is typically the case with solvable groups. On the other hand the Lie rank is sometimes necessary in 
classification problems: it leads to the notions of thin/quasi-thin/generic groups which are essential in the 
Classification of the Finite Simple Groups, as well as in algebraic groups. 

Now the only quasi-simple algebraic groups of Zariski dimension $3$ are of the form $\PSL_2$ or $\SL_2$; 
they also are the only quasi-simple algebraic groups of Lie rank $1$. On the other hand $\PSL_2$ 
and $\SL_2$ are the only ``small'' quasi-simple algebraic groups from a {\em purely group-theoretic} point of view, 
namely a ``local solvability'' condition: the normalizer of each infinite solvable subgroup remains solvable. 
Hence an algebraic group with the latter property must have small Lie rank. 
In the present paper we give a precise meaning to, and prove, such a statement in a much more general context.

Our framework will be that of groups of finite Morley rank, for the global theory of which we refer to 
\cite{BorovikNesin(Book)94} or \cite{AltBorCher(Book)}. Briefly put, groups of finite Morley rank are groups 
equipped with a rudimentary notion of dimension on their first-order definable 
subsets, called the Morley rank for historical reasons in model theory. Since the Morley rank 
satisfies basic axioms reminiscent of the Zariski dimension of algebraic varieties over 
algebraically closed fields, groups of finite Morley rank generalize algebraic groups 
over algebraically closed fields. Conversely, a major 
question which has stirred a huge body of work is the Cherlin-Zilber algebraicity conjecture, 
which postulates that infinite simple groups of finite Morley rank are in fact isomorphic to algebraic 
groups over algebraically closed fields. The algebraicity conjecture holds true at least of groups 
containing an infinite elementary abelian $2$-group \cite{AltBorCher(Book)}; 
so far, the proof is part of the Borovik program for groups with involutions, 
and based on ideas modelled on the Classification of the Finite Simple Groups. 
On the other hand there are potential configurations of simple groups of finite Morley without involutions 
for which the Borovik program is helpless. Here we may refer to the configurations of simple 
``bad" groups of Morley rank $3$ discovered in \cite{Cherlin79}, or more generally to the 
``full Frobenius" groups of finite Morley rank studied in \cite{Jaligot01}. 

In the context of groups of finite Morley, the ``local solvability'' condition mentioned above 
is equivalent to the following. 

\begin{definition}
A group of finite Morley rank is {\em $*$-locally$_{\circ}\o$ solvable} if $N\o(A)$ is solvable for 
each nontrivial abelian connected definable subgroup $A$.
\end{definition}

We now relate our group-theoretic notion of smallness to an abstract version 
of the Lie rank as follows. Since there is a priori no satisfactory first-order analogue of the notion 
of an algebraic torus, we shall deal with certain torsion subgroups throughout. 
Although basic matters such as the conjugacy and structure of 
Sylow $p$-subgroups of a group of finite Morley rank are not settled in general, enough is known 
about abelian divisible $p$-subgroups, which are called {\em $p$-tori}. The maximal ones are conjugate \cite{Cherlin05}, 
and direct powers of a \emph{finite} number of copies of the Pr\"ufer $p$-group $\Z_{p^\infty}$ \cite{BorovikPoizat90}. 
The latter number is called the {\em Pr\"ufer $p$-rank} of the ambient group. 
This will, quite naturally, be our analogue of the Lie rank. 
Indeed, it can easily be seen that in a quasi-simple algebraic group the Lie 
rank and the Pr\"ufer $p$-rank agree: maximal tori of $\SL_2$ or $\PSL_2$ are of dimension $1$ and of Pr\"{u}fer $p$-rank $1$ 
for any prime $p$ different from the characteristic of the ground field. 
We note that in our more abstract context we then have a notion of the Lie rank for 
each prime $p$ such that the ambient group contains a non-trivial divisible abelian $p$-subgroup. 

If $S$ is an abelian $p$-group  for some prime $p$ and $n$ is a natural number, then we denote by $\Omega_{n}(S)$ the subgroup of $S$ generated by all elements of order $p^{n}$. In technical terms that we will define shortly, 
our main theorem takes the following form.

\begin{theorem}%\label{TheoPruferRankspStrongEmbed}
Let $G$ be a connected nonsolvable $*$-locally$_{\circ}\o$-solvable group of finite Morley rank
of Pr\"ufer $p$-rank at least $2$ for some prime $p$, and fix a maximal $p$-torus
$S$ of $G$. Assume that every proper definable connected
subgroup containing $S$ is solvable, that elements of $S$ of order $p$ are
not exceptional, and let
$$B=\<C\o(s)~|~s\in {\Omega_{1}(S)\setminus \{1\}} \>.$$
Then:
\begin{itemize}
\item[$(1)$]
either $B<G$, in which case $B$ is a Borel subgroup of $G$; and if in addition $S$ is a Sylow $p$-subgroup of $N_{N(B)}(S)$, $N(B)$ is $p$-strongly embedded in $G$,
% \item[$(1)$]
% either $B<G$, in which case $B$ is a Borel subgroup of $G$, and moreover
% $N(B)$ is $p$-strongly embedded in $G$ assuming additionally that
% $S$ is a Sylow $p$-subgroup of $N_{N(B)}(S)$,
\item[$(2)$]
or $B=G$, in which case $S$, or equivalently $G$, has Pr\"ufer $p$-rank $2$.
\end{itemize}
\end{theorem}

According to the algebraicity conjecture for simple groups of finite Morley rank, or rather a consequence of it, 
the Pr\"{u}fer $p$-rank of a connected nonsolvable $*$-locally$_{\circ}\o$-solvable group of finite Morley rank 
should be $1$. Hence our Main Theorem deals with configurations 
which are not actually known to exist, but the dichotomy it gives severely limits possibilities in 
both cases. It is obvious in the second case, and in the first case it suffices to recall from \cite{DeloroJaligotI} 
that a definable subgroup $M$ of a group $G$ is called {\em $p$-strongly embedded} if it contains 
non-trivial $p$-elements and $M\cap M^g$ contains none for any $g$ in $G\setminus M$. 
This mimics a similar notion in the theory of finite groups which had been crucial with $p=2$ in the 
Classification of the Finite Simple Groups. 
In any case we note that for groups of finite Morley rank, and typically for $p\neq 2$, 
the ``bad" or ``full Frobenius" groups mentioned above fit in case $(1)$ of our Main Theorem. 

The present paper is actually part of a series which aims at classifying configurations of nonsolvable 
$*$-locally$_{\circ}\o$ solvable groups of finite Morley rank. 
For much more details on such groups, including a few historical remarks, 
we refer to the preliminary article \cite{DeloroJaligotI} of our series. 
We simply recall that this classification started in \cite{CherlinJaligot2004} in the case of 
{\em minimal connected simple} groups: these are the infinite simple groups of finite Morley rank 
all of whose proper definable connected subgroups are solvable (as in $\PSL_2$). 
After a series of generalizations of the original classification of \cite{CherlinJaligot2004} it was realized 
that much of the theory of minimal connected simple groups transfers readily to {\em $*$-locally$\o$-solvable} groups, 
which are defined as the $*$-locally$_\circ\o$-solvable ones but where the main condition 
that $N\o(A)$ is solvable is required for {\em any} nontrivial abelian subgroup $A$. 
This shift from minimal connected simple groups to $*$-locally$\o$-solvable groups corresponds, 
in finite group theory, to a shift from the {\em minimal simple} groups studied for the 
Feit-Thompson (Odd Order) Theorem to the {\em $N$-groups} classified later by Thompson. 
Of the two main variations on the notion of ``local solvability" from \cite{DeloroJaligotI} we shall work here with
the most general $*$-local$_\circ\o$-solvability, which allows $\SL_2$ in addition to $\PSL_2$. 
Notice that if $G$ is $*$-locally$_\circ\o$-solvable, then it can contain elements $x$ (of order necessarily finite)
with $C_G\o (x)$ nonsolvable. Following \cite{DeloroJaligotI} such elements $x$ are 
called {\em exceptional}, refering somehow to the central involution of $\SL_2$. 

Since \cite{CherlinJaligot2004} it seems unrealistic to hope for a complete reduction of  
$*$-locally$_\circ\o$ solvable groups to the algebraic ones, even assuming the presence of involutions 
in order to proceed to a much sharper analysis. However a reduction to a very small number of configurations 
will be obtained in \cite{DeloroJaligotII}, along the lines of \cite{CherlinJaligot2004} and subsequent papers. 
This can be seen as part of the Borovik program for classifying simple groups with involutions, in the utterly 
critical case of minimal configurations. 
We simply note that the present part of the analysis does not depend on $p$ and can be reached by very 
general means (and some of them, such as signalizer functors, could be borrowed from finite group theory).

We also note that our Main Theorem generalizes the dichotomy represented by Sections 6 and 7 of 
\cite{CherlinJaligot2004}, concerning $p=2$ in minimal connected simple groups also satisfying 
a simplifying ``tameness" assumption. Sections 6 (resp. 7) there corresponds, in Pr\"ufer $2$-rank at least $2$, 
to $C_G\o (\Omega_1 (S))$ not being (resp. being) a Borel subgroup, two cases corresponding respectively 
to our cases $(2)$ and $(1)$ here. Then came \cite{BurdgesCherlinJaligot07} where our corresponding case $(1)$ 
was shown not to exist, still for minimal connected simple groups and for $p=2$, 
but without the ``tameness" assumption. In \cite{DeloroJaligotII} we will reach essentially the same conclusions in the 
much more general context of $*$-locally$_{\circ}\o$-solvable groups, applying the general dichotomy of 
the present paper with $p=2$. In particular we will bound by $2$, in full generality, the Pr\"ufer $2$-rank of a 
nonsolvable $*$-locally$_\circ\o$ solvable group. 
But in any case, case $(2)$ of our Main Theorem still stands around, 
even with $p=2$ and in the context of tame minimal connected simple groups of \cite{CherlinJaligot2004} where the configuration 
is described with high precision. 

We collect some raw material in \S\ref{S:Preliminaries} and the proof of our Main Theorem takes place in \S\ref{S:Proof}. In \S\ref{S:Final} we will make more comments on the difficulty to deal with the configuration arising in case $(2)$; we will also give a form of our Main Theorem more directly applicable in
\cite{DeloroJaligotII} (essentially explaining how to deal which the two extra assumptions that every proper definable connected subgroup containing $S$ is solvable and that elements of $S$ of order $p$ are
not exceptional).

\section{Preliminaries}\label{S:Preliminaries}

For general reference on groups of finite Morley rank we refer to \cite{BorovikNesin(Book)94} or \cite{AltBorCher(Book)}, and for more specific facts about $*$-locally$_\circ\o$ solvable groups of finite Morley rank we refer to \cite{DeloroJaligotI}. We simply recall that groups of finite Morley rank satisfy the descending chain condition on definable subgroups, and that any group $G$ of finite Morley rank has a {\em connected component}, i.e., a smallest (normal) definable subgroup of finite index, denoted by $G\o$.

\subsection{Unipotence theory}

A delicate point in our proof is the use of the abstract unipotence theory for groups of finite Morley rank. We follow the general treatment of \cite{FreconJaligot07} and \cite[\S 2.1]{DeloroJaligotI}, and for the sake of self-containment we shall try and put unipotence theory in a nutshell. $\mathcal{P}$ denotes the set of prime numbers.

\begin{definition}
A {\em unipotence parameter} is a pair $\tilde{p} = (p, r) \in (\{\infty\} \cup \mathcal{P}) \times (\N \cup \{\infty\})$ with $p < \infty$ if and only if $r = \infty$. In addition, $p$ is called the {\em characteristic} of $\tilde{p}$ and $r$ is called the {\em unipotence degree} of $\tilde{p}$.
\end{definition}

Such a notion enables a parallel treatment to two theories which had been considered distinct, namely $p$-unipotence theory ($p$ a prime), and Burdges' more recent notion of null characteristic graded unipotence which originated in \cite{Burdges03}. For $p$ a prime, a {\em $p$-unipotent} subgroup is by definition a nilpotent definable connected $p$-group of bounded exponent, and for $p=\infty$ a definition is given in terms of generation by certain definable connected abelian subgroups with torsion-free quotients of rank $r$
\cite[\S2.3]{FreconJaligot07}. In the extreme case where $(p,r) = (\infty,0)$ one makes use of so-called decent tori \cite{Cherlin05}.
All this gives rise, for any group $G$ of finite Morley rank, to the group $U_{\tilde p}(G)$ as the group generated by subgroups as above. This group is always definable and connected, and for $p$ a prime $U_{(p,\infty)}(G)$ is simply denoted by $U_p(G)$.

Imposing nilpotence, it leads to a general notion of $\tilde p$-subgroup for any unipotence parameter $\tilde p$, and now imposing maximality with respect to inclusion, it provides a notion of {\em Sylow $\tilde p$-subgroup} \cite[\S2.4]{FreconJaligot07}, with properties much similar to those of Sylow $p$-subgroups of finite groups.

In groups of finite Morley rank, the {\em Fitting} subgroup, i.e., the group generated by all normal nilpotent subgroups, is always definable and nilpotent, and in particular the {\em unique} maximal normal nilpotent subgroup. If one lets $d_p (H)$ be the maximal $r$ such that $H$ contains a non-trivial $(p, r)$-unipotent subgroup \cite[Definition 2.5]{DeloroJaligotI}, then the interest for the unipotence theory comes from the following in connected solvable groups.

\begin{fact}\label{FactUnipHeaviest}
{\bf \cite[Fact 2.8]{DeloroJaligotI}}
Let $H$ be a connected solvable group of finite Morley rank and $\tilde{p}=(p,r)$
a unipotence parameter with $r>0$. Assume $d_{p}(H)\leq r$. Then
$U_{\tilde p}(H)\leq F\o(H)$, and in particular $U_{\tilde p}(H)$ is nilpotent.
\end{fact}

\subsection{Uniqueness Theorem in $*$-locally$_{\circ}\o$-solvable groups}

The main result about $*$-locally$_{\circ}\o$-solvable groups is a Uniqueness Theorem analogous to a similar result of Bender about minimal simple groups in finite group theory. We give here the most general statement as proved in \cite{DeloroJaligotI}. We will not use it as such, but the proof of our Main Theorem will be based on an elaborate version of this Uniqueness Theorem.

\begin{fact}\label{UniquenessLemma}
{\bf \cite[Uniqueness Theorem 4.1]{DeloroJaligotI}}
Let $G$ be a $*$-locally$_{\circ}\o$-solvable group of finite Morley rank,
$\tilde p=(p,r)$ a unipotence parameter with $r>0$,
and $U$ a Sylow $\tilde p$-subgroup of $G$.
Assume that $U_{1}$ is a nontrivial definable $\tilde p$-subgroup of $U$ containing
a nonempty (possibly trivial) subset $X$ of $G$ such that $d_{p}(C\o(X))\leq r$.
Then $U$ is the unique Sylow $\tilde p$-subgroup of $G$
containing $U_{1}$, and in particular $N(U_{1})\leq N(U)$.
\end{fact}

A {\em Borel} subgroup is a maximal definable connected solvable subgroup. Fact \ref{UniquenessLemma} has the following consequence on Borel subgroups of $*$-locally$_{\circ}\o$-solvable groups of finite Morley rank.

\begin{fact}\label{corUnicityBorelSuperMax}
{\bf \cite[Corollary 4.4]{DeloroJaligotI}}
Let $G$ be a $*$-locally$_{\circ}\o$-solvable group of finite Morley rank,
$\tilde p=(p,r)$ a unipotence parameter with $r>0$ such that $d_{p}(G)=r$.
Let $B$ be a Borel subgroup of $G$ such that $d_{p}(B)=r$. Then
$U_{\tilde p}(B)$ is a Sylow $\tilde p$-subgroup of $G$,
and if $U_{1}$ is a nontrivial definable $\tilde p$-subgroup of $B$, then
$U_{\tilde p}(B)$ is the unique Sylow $\tilde p$-subgroup of $G$
containing $U_{1}$, $N(U_{1})\leq N(U_{\tilde p}(B))=N(B)$, and $B$ is the unique
Borel subgroup of $G$ containing $U_{1}$.
\end{fact}

\subsection{Torsion}

The results we shall use about torsion are rather elementary. We first give the general decomposition of nilpotent groups of finite Morley rank, in terms of both classical Sylow $p$-subgroups and Sylow $\tilde p$-subgroups.

\begin{fact}\label{StrucNilpGroups}
{\bf \cite[Fact 2.3]{DeloroJaligotI}}
Let $G$ be a nilpotent group of finite Morley rank.
\begin{itemize}
\item[$(1)$]
$G$ is the central product of its
Sylow $p$-subgroups and its Sylow $(\infty,r)$-subgroups, which are all divisible for $\tilde{p}$ of the form $(\infty,r)$ \cite[p. 16]{FreconJaligot07}.
\item[$(2)$]
If $G$ is connected, then $G$ is the central
product of its Sylow $\tilde p$-subgroups.
\end{itemize}
\end{fact}

Throughout the rest of this subsection, $p$ denotes a prime number. For an arbitrary subgroup $S$ of a group of finite Morley rank, one defines the (generalized) {\em connected component} of $S$ by $S\o=H\o(S)\cap S$, where $H(S)$ is the smallest definable subgroup of $G$ containing $S$. (Of course, $H(S)$ exists by descending chain condition on definable subgroups.) It is easily checked that $S\o$ has finite index in $S$.

\begin{fact}\label{StructpSylSolvGps}
{\bf \cite[Corollary 6.20]{BorovikNesin(Book)94}}
Let $p$ be a prime and $S$ a $p$-subgroup of a solvable group of finite Morley rank,
or more generally a locally finite $p$-subgroup of any group of finite Morley rank. Then:
\begin{itemize}
\item[$(1)$]
$S\o$ is a central product of a $p$-torus and a $p$-unipotent subgroup.
\item[$(2)$]
If $S$ is infinite and has bounded exponent, then $Z(S)$ contains infinitely many
elements of order $p$.
\end{itemize}
\end{fact}

We also recall that Sylow $p$-subgroups are conjugate in solvable groups of finite Morley rank
\cite[Theorem 9.35]{BorovikNesin(Book)94}. Furthermore they are connected in connected solvable groups of finite Morley rank by \cite[Theorem 9.29]{BorovikNesin(Book)94}; in this case a Sylow $p$-subgroup $S$ satisfies $S=S\o$ and has a decomposition as in Fact \ref{StructpSylSolvGps} $(1)$, a point that will be frequently used below.

If $G$ is a group of finite Morley rank, we denote by $O_{p'}(H)$
the largest normal definable connected subgroup without $p$-torsion. It exists by ascending
chain condition on definable connected subgroups and the following elementary property of lifting
of torsion valid in groups of finite Morley rank (or more generally in groups with the descending chain condition on definable subgroups):

\begin{fact}\label{LiftingTorsion}
{\bf \cite[Lemma 2.18]{AltBorCher(Book)}}
Let $G$ be a group of finite Morley rank,
$N$ a normal definable subgroup of $G$, and $x$ an element of $G$ such that $x$
has finite order $n$ modulo $N$. Then the coset $xN$
contains an element of finite order, involving the same prime divisors as  $n$.
\end{fact}

The following will be useful when dealing with $p$-strongly embedded subgroups.

\begin{lem}\label{FaitH/Op'HDivAbelian}
{\bf (Compare with \cite[Lemma 3.2]{CherlinJaligot2004})}
Let $H$ be a connected solvable group of finite Morley rank
such that $U_{p}(H)=1$. Then $H/O_{p'}(H)$ is divisible abelian.
\end{lem}
\proof
Dividing by $O_{p'}(H)$, we may assume it is trivial and we want
to show that $H$ is divisible abelian.

Let $F=F\o(H)$. As $O_{p'}(H)=1$, $O_{p'}(F)=1$ as well,
and $U_{q}(H)=1$ for any prime $q$ different from $p$. By assumption
$U_{p}(H)=1$ also, and $F$ is divisible by Fact \ref{StrucNilpGroups}.
As $F'$ is torsion-free, by \cite[Theorem 2.9]{BorovikNesin(Book)94}
or Fact \ref{StrucNilpGroups} and \cite[Corollary 2.2]{DeloroJaligotI},
it must be trivial by assumption. (In connected groups of finite Morley rank, derived subgroups are definable and connected by a well-known corollary of Zilber's indecomposability theorem.)
Hence $F$ is divisible abelian.

To conclude it suffices to show that $F$ is central in $H$, as then $H$
is nilpotent, hence equal to $F$, and hence divisible abelian, as desired.
Let $h$ be any element of $H$; we want to show that $[h,F]=1$. Notice that that map
$f \mapsto [h,f]$, from $F$ to $F$, is a definable group homomorphism.
As the torsion subgroup of $F$ is central in $H$ (by \cite[Fact 2.7 $(1)$]{DeloroJaligotI}, or using
\cite[Fact 2.1]{DeloroJaligotI}), it is contained in the kernel of the previous map and
Fact \ref{LiftingTorsion} shows that the image of the previous map, i.e., $[h,F]$, is torsion-free.
Hence $[h,F]\leq O_{p'}(F)=1$, as desired.
\qed

\subsection{Generation by centralizers}

In the present subsection we prove miscellaneous lemmas concerning generation by centralizers.

\begin{lem}\label{FaitDefGpAutpGp}
Let $\tilde p$ be a unipotence parameter and $q$ a prime number.
Let $H$ be a $\tilde p$-group of finite Morley rank without elements of
order $q$, and assume $K$ is a definable solvable $q$-group of automorphisms of $H$
of bounded exponent. Then $C_{H}(K)$ is a definable $\tilde p$-subgroup of $H$.
\end{lem}
\proof
By descending chain condition on centralizers, $C_{H}(K)$ is the centralizer of a
finitely generated subgroup of $K$, and by local finiteness of the latter we may assume
$K$ finite. In particular $C_H(K)$ is connected by \cite[Fact 3.4]{Burdges03}.

When $\tilde p=(\infty,0)$, $H$ is a good torus in the sense of \cite{Cherlin05}, and in particular
$(\infty,0)$-homogeneous in the sense of \cite[Lemma 2.17]{FreconJaligot07}, and the connected subgroup $C_H(K)$ is
also a good torus. Otherwise, $C_H(K)$ is also a $\tilde p$-group, by \cite[Lemma 3.6]{Burdges03} when the
unipotence parameter is finite, or by \cite[Lemma 2.17-c]{FreconJaligot07} when the
characteristic is finite.
\qed

\begin{fact}\label{BiGenerationWithpElts}
{\bf \cite[Fact 3.7]{Burdges03}}
Let $H$ be a solvable group of finite Morley rank without elements of order $p$
for some prime $p$.
Let $E$ be a finite elementary abelian $p$-group acting definably on $H$.
Then:
$$H=\<C_{H}(E_{0})~|~E_{0}\leq E,~[E:E_{0}]=p\>.$$
\end{fact}

We recall that a {\em Carter} subgroup of a group of finite Morley is, by definition, a definable connected
nilpotent subgroup of finite index in its normalizer.

\begin{lem}\label{Bipgeneration}
Let $H$ be a connected solvable group of finite Morley rank
such that $U_{p}(H)=1$ for some prime $p$. Suppose that $H$ contains an elementary abelian
$p$-group $E$ of order $p^{2}$. Then:
$$H=\<C\o_{H}(E_{0})~|~E_{0}
\mbox{~is a cyclic subgroup of order $p$ of~}E\>.$$
\end{lem}
\proof
By assumption, Fact \ref{StructpSylSolvGps}, and \cite[Theorem 9.29]{BorovikNesin(Book)94},
Sylow $p$-subgroups of $H$ are $p$-tori. Hence $E$ is in a maximal $p$-torus of $H$,
which is included in a Carter subgroup $Q$ of $H$ by
\cite[Theorem 3.3]{FreconJaligot07}.
By Lemma \ref{FaitH/Op'HDivAbelian},
$H/O_{p'}(H)$ is abelian. As Carter subgroups cover all abelian quotients
in connected solvable groups of finite Morley rank by \cite[Corollary 3.13]{FreconJaligot07},
$H=O_{p'}(H)\cdot Q$. As $E\leq Z(Q)$, it suffices to show that:
$$O_{p'}(H)=\<C\o_{O_{p'}(H)}(E_{0})~|~E_{0}
\mbox{~is a cyclic subgroup of order $p$ of~}E\>.$$
But the generation by the full centralizers is given by
Fact \ref{BiGenerationWithpElts}, and these centralizers are connected
by \cite[Fact 3.4]{Burdges03}.
\qed

\bigskip
A subgroup is called {\em $p$-toral}, or just {\em toral}, if it is contained in a $p$-torus of the ambient group.

\begin{lem}\label{CorBipgeneration}
Let $H$ be a connected solvable group of finite Morley rank with an elementary abelian $p$-toral
subgroup $E$ of order $p^2$ for some prime $p$. Then:
$$H=\<C\o_{H}(E_{0})~|~E_{0}
\mbox{~is a cyclic subgroup of order $p$ of~}E\>.$$
\end{lem}
\proof
For a connected nilpotent group of finite Morley rank $L$, we define
the ``complement" $C_{p}(L)$ of $U_{p}(L)$, namely the product
of all factors of $L$ as in Fact \ref{StrucNilpGroups} $(2)$, except $U_p(L)$.

Now if $H$ is any connected solvable group of finite Morley rank and
$Q$ a Carter subgroup of $H$, then $H=QF\o(H)$ by \cite[Corollaire 3.13]{FreconJaligot07}, and $H$ is the product of the definable
connected subgroup $C_{p}(Q)C_{p}(F\o(H))$ with the normal definable
connected subgroup $U_{p}(H)$, and the first factor has trivial
$p$-unipotent subgroups.

In our particular case, $E$ is by torality contained in a $p$-torus, and the latter
is contained in a Carter subgroup $Q$ of $H$ by \cite[Th\'eor\`eme 3.3]{FreconJaligot07}. By
\cite[Theorem 9.29]{BorovikNesin(Book)94} and Fact \ref{StructpSylSolvGps},
$E$ centralizes the normal definable connected subgroup $U_{p}(H)$,
so it suffices to show the generation by the connected components of centralizers in
$C_{p}(Q)C_{p}(F\o(H))$. But this follows from
Lemma \ref{Bipgeneration}.
\qed

%\subsection{Pr\"ufer ranks and strong embedding}\label{SectionPruferRanks}

\section{Proof of our theorem}\label{S:Proof}

We now turn to proving our Main Theorem which we restate.

\begin{theorem}%\label{TheoPruferRankspStrongEmbed}
Let $G$ be a connected nonsolvable $*$-locally$_{\circ}\o$-solvable group of finite Morley rank
of Pr\"ufer $p$-rank at least $2$ for some prime $p$, and fix a maximal $p$-torus
$S$ of $G$. Assume that every proper definable connected
subgroup containing $S$ is solvable, that elements of $S$ of order $p$ are
not exceptional, and let
$$B=\<C\o(s)~|~s\in {\Omega_{1}(S)\setminus \{1\}} \>.$$
Then:
\begin{itemize}
\item[$(1)$]
either $B<G$, in which case $B$ is a Borel subgroup of $G$; and if in addition $S$ is a Sylow $p$-subgroup of $N_{N(B)}(S)$, $N(B)$ is $p$-strongly embedded in $G$,
% \item[$(1)$]
% either $B<G$, in which case $B$ is a Borel subgroup of $G$, and moreover
% $N(B)$ is $p$-strongly embedded in $G$ assuming additionally that
% $S$ is a Sylow $p$-subgroup of $N_{N(B)}(S)$,
\item[$(2)$]
or $B=G$, in which case $S$ has Pr\"ufer $p$-rank $2$.
\end{itemize}
\end{theorem}

\paragraph{Here begins the proof.}
Let $M=N(B)$. As $B$ is definable (and connected) by Zilber's generation lemma, $M$ is definable as well.
As $B$ contains a generous Carter subgroup $Q$ of $G$ containing $S$ by
\cite[Fact 3.32]{DeloroJaligotI}, the conjugacy of generous Carter subgroups of
\cite{Jaligot06} and a Frattini argument give
$M=N(B)\subseteq BN(Q)$, and as $Q$ is almost selfnormalizing
$B=M\o$. (With the same notation, this holds of course for an arbitrary $p$-torus $S$
in an arbitrary group $G$ of finite Morley rank.)

\subsection{Getting rid of $(1)$}

Assume first
$$B<G.\leqno(1)$$

By assumption $B\leq B_{1}$ for some Borel subgroup $B_{1}$ of $G$.
As $S\leq B \leq B_1$, Lemma \ref{CorBipgeneration}
implies that $B=B_{1}$, and thus $B$ is a Borel subgroup of $G$.

In particular, Sylow $p$-subgroups of $B$ are conjugate, as in any solvable group of finite Morley rank.

We now make the new assumption that $S$ is a Sylow $p$-subgroup of
$N_{N(B)}(S)$.

\begin{lem}\label{Lem3.1DJI'}
$U_{p}(C(s))=1$ for every element $s$ of order $p$ of $S$.
\end{lem}
\proof
It suffices to apply Fact \ref{StructpSylSolvGps} in the connected solvable group $B$, again with the fact that Sylow $p$-subgroups of connected solvable groups of finite Morley rank are connected.
\qed

\bigskip
We claim that $M=N(B)$ is $p$-strongly embedded in $G$ in this case by using
a ``black hole" principle (a term going back to Harada) similar to the one
used in \cite[\S2.2]{BurdgesCherlinJaligot07}, and already contained in
\cite[Lemma 7.3]{CherlinJaligot2004}.
We note that Lemma \ref{Lem3.1DJI'} implies that $S$ is a
Sylow $p$-subgroup of $B$ indeed, and of
$M$ as well, as $M=BN_{N(B)}(S)$ by a Frattini Argument. In particular
$M/B$ has trivial Sylow $p$-subgroups by lifting of torsion, Fact \ref{LiftingTorsion}.

Assume that $M\cap M^{g}$ contains an element $s$ of order
$p$ for some $g$ in $G$. Notice that $s$ is actually in $B\cap B^g$,
and $p$-toral. By connectedness and conjugacy of Sylow $p$-subgroups in
connected solvable groups, the definition of $B$ implies that $C\o(s')\leq B$
for any element $s'$ of order $p$ of $B$. Similarly, $C\o(s')\leq B^{g}$ whenever
$s'$ has order $p$ and is in $B^g$. By conjugacy in $B$ we may assume $s$ in $S$,
and $\Omega_1(S)\leq S\leq C\o(s)\leq B\cap B^g$.
By Lemma \ref{Bipgeneration} or \ref{CorBipgeneration}
applied in $B$ and in $B^{g}$ we get
$B^{g}={\<C\o(s)~|~s\in \Omega_{1}(S)\setminus \{1\} \>}=B$.
Thus $g$ normalizes $B$, and is in $M$.

Hence $M=N(B)$ is $p$-strongly embedded in $G$ under the extra assumption adopted here,
and this proves clause $(1)$ of the Main Theorem. %\ref{TheoPruferRankspStrongEmbed}

\subsection{Case $(2)$; $p$-elements}

We now pass to the second case
$$B=G.\leqno(2)$$
We will eventually show that clause $(2)$ of
our Main Theorem %\ref{TheoPruferRankspStrongEmbed}
holds by reworking the begining
of Section 6 of \cite{CherlinJaligot2004}. We first put aside $p$-unipotent subgroups.

\begin{lem}\label{LemAnyBorelContspPerp}
Any Borel subgroup containing a toral element of order $p$ has trivial
$p$-unipotent subgroups.
\end{lem}
\proof
Assuming the contrary, we may assume after conjugacy of decent tori \cite{Cherlin05} that a
Borel subgroup $L$ with $U_{p}(L)$ nontrivial contains an element $s$ of $S$
of order $p$. Then $U_{p}(C(s))$ is nontrivial by Fact \ref{StructpSylSolvGps},
contained in a unique Borel subgroup $B_{1}$ of $G$ by Fact \ref{corUnicityBorelSuperMax}.
(Actually $B_{1}=L$.) By
Fact \ref{corUnicityBorelSuperMax}, $B_{1}$ is the unique Borel subgroup
containing any given nontrivial $p$-unipotent subgroup of $U_{p}(C(s))$.
Now any element $s'$ of order $p$ of $S$ normalizes $U_{p}(C(s))$, and thus
$U_{p}(C(s,s'))\neq 1$ by
Fact \ref{StructpSylSolvGps}, and as $B_{1}$ is the unique Borel subgroup
containing the latter group we get $C\o(s')\leq B_{1}$. This shows that $B\leq B_{1}$,
a contradiction as $B=G$ is nonsolvable under the current assumption.
\qed

\bigskip
In other words, nontrivial $p$-toral elements commute with no nontrivial $p$-unipotent
subgroups. This can be stated more carefully as follows.

\begin{cor}\label{CorLemAnyBorelContspPerp}
Any connected solvable subgroup which is $\<s\>$-invariant for some $p$-toral element
$s$ of order $p$ has trivial $p$-unipotent subgroups.
\end{cor}
\proof
Otherwise $s$ would normalize a nontrivial $p$-unipotent
subgroup, and by Fact \ref{StructpSylSolvGps} it would
centralize a nontrivial $p$-unipotent subgroup.
\qed

\bigskip
Our assumption $(2)$ on $B$ yields similarly a property antisymmetric to
the black hole principle implied by assumption $(1)$.
Let $E$ denote the elementary abelian $p$-group $\Omega_1(S)$.

\begin{lem}\label{LemDispatchCentInvBorels}
Let $E_1$ be a subgroup of $E$ of order at least $p^{2}$. Then for any proper
definable connected
subgroup $L$ there exists an element $s$ of order $p$ of $E_{1}$ such that
$C\o(s)\nleq L$. In particular $G=\<C\o(s)~|~s\in E_{1}\setminus \{1\}\>$.
\end{lem}
\proof
Assume on the contrary $C\o(s)\leq L$ for any element $s$ of
order $p$ of $E_1$.

We claim that $C\o(t)\leq L$ for any element $t$ of order $p$ of $E$.
In fact, as $E_1\leq S\leq C\o(t)$, $C\o(t)$ is by Lemma \ref{Bipgeneration}
generated by its subgroups of the form $C\o(t,s)$, with
$s$ of order $p$ in $E_1$. As these groups are all contained in $L$ by assumption,
our claim follows.

Hence we have $B\leq L<G$. But under our current assumption
$B=G$, and this is a contradiction.

Our last claim follows immediately.
\qed

\begin{cor}\label{CorCs<Cs}
There exists an element $s$ of order $p$ of $E$ such that $C\o(S)<C\o(s)$.
\end{cor}
\proof
$C\o(S)$ is $S$-local$\o$, and thus solvable by $*$-local$_{\circ}\o$ solvability of $G$
and \cite[Lemma 3.4]{DeloroJaligotI}. As $C\o(S)\leq C\o(s)$ for
any element $s$ of order $p$ of $S$, it suffices to apply
Lemma \ref{LemDispatchCentInvBorels}.
\qed

\bigskip
Recall from \cite[Definition 2.5]{DeloroJaligotI} that for any group $G$ of finite Morley rank, $d(G)$ denotes the maximum of $d_{\infty}(G)$ and of $\max_{p\in {\cal P}}(d_p(G))$, i.e., the infinite symbol if $G$ contains a nontrivial $p$-unipotent subgroup, and otherwise the maximal $r$ in $\N$ such that $G$ contains a nontrivial $(\infty,r)$-subgroup. This number exists and belongs to $\N \sqcup{\{\infty\}}$ when $G$ is infinite \cite[Lemma 2.6]{DeloroJaligotI}.

\begin{lem}\label{LemExistsLdOp'L1}
There exists an element $s$ of order $p$ of $E$ such that
$$d(O_{p'}(C\o(s)))\geq 1.$$
\end{lem}
\proof
Assume the contrary, and let $s$ be an arbitrary element of order $p$ of $E$.
By our assumption that $d(O_{p'}(C\o(s)))\leq 0$,
$O_{p'}(C\o(s))$ is trivial or a good torus by
\cite[Lemma 2.6]{DeloroJaligotI}, and
central in $C\o(s)$ by \cite[Fact 2.7 $(1)$]{DeloroJaligotI}. Notice that
$U_{p}(C\o(s))=1$ by Lemma \ref{LemAnyBorelContspPerp}. As
$C\o(s)/{O_{p'}(C\o(s))}$ is abelian by Lemma \ref{FaitH/Op'HDivAbelian},
$C\o(s)$ is nilpotent. Now $S$ is central in $C\o(s)$ by Fact \ref{StrucNilpGroups} $(2)$.
In particular $C\o(S)=C\o(s)$, and this holds for any element $s$ of order $p$
of $E$. We get a contradiction to Corollary \ref{CorCs<Cs}.
\qed

\bigskip
It follows in particular from Lemma \ref{LemExistsLdOp'L1}
that there exist definable connected subgroups $L$
containing $C\o(s)$ for some element $s$ of order $p$ of $E$
and such that $O_{p'}(L)$ is not a good torus.
Choose then a unipotence parameter $\tilde q=(q,r)$
different from $(\infty,0)$ such that $r$ is maximal in the set of all
$d_{q}(O_{p'}(L))$, where $L$ varies in the set of all definable connected
solvable subgroups with the above property.

Notice that there might exist several such maximal unipotence parameters
$\tilde q$, maybe one for $q=\infty$ and several ones for $q$ prime, except
for $q=p$ by Corollary \ref{CorLemAnyBorelContspPerp}.

It will also be shortly and clearly visible below that the notion of
maximality for $\tilde q$ is the same when $L$ varies in two smaller subsets
of all definable connected solvable subgroups
containing $C\o(s)$ for some $s$ of order $p$ of $E$: the set of
Borel subgroups with this property
on the one hand, and exactly the finite set of subgroups of the form $C\o(s)$ on the other.

\begin{lem}
Let $L$ be any definable connected solvable subgroup containing $C\o(s)$
for some element $s$ of order $p$ of $E$. Then $U_{\tilde q}(O_{p'}(L))$
is a normal definable connected nilpotent subgroup of $L$.
\end{lem}
\proof
As $O_{p'}(L)$ is normal in $L$, it suffices to show that its definably
characteristic subgroup $U_{\tilde q}(O_{p'}(L))$ is nilpotent. But the latter is in
$F(O_{p'}(L))$ by Fact \ref{FactUnipHeaviest} and the maximality of $r$.
\qed

\begin{cor}\label{CorEveryThingInFOp'}
Let $L$ be any definable connected solvable subgroup containing $C\o(s)$
for some element $s$ of order $p$ of $E$. Then any definable $\tilde q$-subgroup
of $L$ without elements of order $p$ is in $U_{\tilde q}(F(O_{p'}(L)))$.
\end{cor}
\proof
Let $U$ be such a subgroup. As $U_{p}(L)=1$ by Lemma \ref{LemAnyBorelContspPerp},
$L/O_{p'}(L)$ is (divisible) abelian by Lemma \ref{FaitH/Op'HDivAbelian}, and
thus $U\leq O_{p'}(L)$, and $U\leq U_{\tilde q}(O_{p'}(L))$. Now it
suffices to apply the normality and the nilpotence of the latter.
\qed

\subsection{A Uniqueness Theorem via elementary $p$-groups}

We now prove a version of the Uniqueness Theorem (Fact \ref{UniquenessLemma})
with a combined action, more precisely where the assumption on unipotence
degrees of centralizers is replaced by an assumption of invariance by
a sufficiently ``large" $p$-toral
subgroup. For this purpose we first note the following.

\begin{lem}\label{LemActionUnipDegOnpPerpByE1}
Let $E_{1}$ be a subgroup of order at least $p^{2}$ of $E$, and $H$ a
definable connected solvable $E_{1}$-invariant subgroup. Then
$d_{q}(O_{p'}(H))\leq r$.
\end{lem}
\proof
Assume toward a contradiction $r'>r$, where $r'$ denotes
$d_{q}(O_{p'}(H))$. In this case $r$ is necessarily finite, and $q=\infty$.
By Fact \ref{FactUnipHeaviest}, $U_{(\infty,r')}(O_{p'}(H))\leq F\o(O_{p'}(H))$,
and this nontrivial definable $(\infty,r')$-subgroup is $E_{1}$-invariant.
Fact \ref{BiGenerationWithpElts} gives an element $s$ of order $p$
in $E_{1}$ such that
$$C_{U_{(\infty,r')}(O_{p'}(H))}(s)\neq 1.$$
But the latter is an $(\infty,r')$-group by
Lemma \ref{FaitDefGpAutpGp}.
Now considering the definable connected solvable subgroup
$C\o(s)$ gives a contradiction to the maximality of $r$, as $C\o(s)/O_{p'}(C\o(s))$
is (divisible) abelian as usual and the centralizer above is
connected without elements of order $p$, and thus contained in $O_{p'}(C\o(s))$.
\qed

\bigskip
As already mentioned around the definition of maximal parameters $\tilde q$ (after Lemma \ref{LemExistsLdOp'L1}),
the same argument shows that $r$ is also exactly
the maximum of the $d_{q}(O_{p'}(L))$ different from $0$, with $L$ varying
in the set of {\em Borel} subgroups containing $C\o(s)$ for some element $s$ of
order $p$ of $E$
(instead of all definable connected solvable subgroups $L$ with the same property),
and similarly with $L$ varying in the set of subgroups $C\o(s)$ for some element $s$ of
order $p$ of $E$.

We now prove our version of the
Uniqueness Theorem, Fact \ref{UniquenessLemma}, specific to the configuration considered here.

\begin{theo}\label{RelUniqLem}
Let $E_{1}$ be a subgroup of order at least $p^{2}$ of $E$. Then any
$E_{1}$-invariant nontrivial definable $\tilde q$-subgroup without elements
of order $p$ is contained in a unique maximal such.
\end{theo}
\proof
Let $U_{1}$ be the $\tilde q$-subgroup under consideration.
Let $U$ be a maximal $E_{1}$-invariant definable $\tilde q$-subgroup
without elements of order $p$ containing $U_{1}$.

Assume $V$ is another such subgroup, distinct from $U$,
and chosen so as to maximize the rank of
$U_{2}=U_{\tilde q}(U\cap V)$. As $1< U_{1}\leq U_{2}$,
the subgroup $U_{2}$ is nontrivial. As $U_{2}$ is nilpotent,
$N:=N\o(U_{2})$ is solvable by $*$-local$_{\circ}\o$ solvability of $G$.
Note that $U_{2}<U$, as otherwise $U=U_{2}\leq V$ and $U=V$ by maximality
of $U$. Similarly $U_{2}<V$, as otherwise $V=U_{2}\leq U$ and $V=U$
by maximality of $V$. In particular, by normalizer condition
\cite[Proposition 2.8]{FreconJaligot07}, $U_{2}<U_{\tilde q}(N_{U}(U_{2}))$ and
$U_{2}<U_{\tilde q}(N_{V}(U_{2}))$.

We claim that $d_{q}(O_{p'}(N))=r$. Actually
$d_{q}(O_{p'}(N))\leq r$ by Lemma \ref{LemActionUnipDegOnpPerpByE1},
and as $O_{p'}(N)$ contains $U_{2}$ which is nontrivial and
of unipotence degree $r$ in characteristic $q$ we get $d_{q}(O_{p'}(N))=r$.

By Fact \ref{FactUnipHeaviest} and the fact that $r\geq 1$
we get $U_{\tilde q}(O_{p'}(N))\leq F\o(O_{p'}(N))$.
In particular $U_{\tilde q}(O_{p'}(N))$ is
nilpotent, and contained in a maximal definable $E_{1}$-invariant
$\tilde q$-subgroup without elements of order $p$, say $\Gamma$.
Notice that $N$, being $E_{1}$-invariant, satisfies $U_{p}(N)=1$, and
$N/O_{p'}(N)$ is abelian as usual.
Now $U_{1}\leq U_{2}<U_{\tilde q}(N_{U}(U_{2}))\leq \Gamma$, so our
maximality assumption implies that $\Gamma=U$. In particular
$U_{\tilde q}(N_{V}(U_{2}))\leq \Gamma =U$. But then
$U_{2}<U_{\tilde q}(N_{V}(U_{2}))\leq U_{\tilde q}(U\cap V)=U_{2}$,
a contradiction.
\qed

\begin{cor}\label{CorUnipTheoAction}
Let $E_{1}$ be a subgroup of order at least $p^{2}$ of $E$.
\begin{itemize}
\item[$(1)$]
If $U_{1}$ is a nontrivial $E_{1}$-invariant definable $\tilde q$-subgroup without
elements of order $p$, then $U_{1}$ is contained in a unique maximal
$E_{1}$-invariant definable connected solvable subgroup $B$.
Furthermore $U_{\tilde q}(O_{p'}(B))$ is the unique
maximal $E_{1}$-invariant definable $\tilde q$-subgroup without elements
of order $p$ containing $U_{1}$, and, for any element $s$ of order $p$
of $E_{1}$ with a nontrivial centralizer in $U_{1}$, $C\o(s)\leq B$ and $B$ is
a Borel subgroup of $G$.
\item[$(2)$]
$U_{\tilde q}(O_{p'}(C\o(E_{1})))$ is trivial.
\end{itemize}
\end{cor}
\proof
$(1)$. Assume $B_{1}$ and $B_{2}$ are two maximal $E_{1}$-invariant definable
connected solvable subgroups
containing $U_{1}$. We have $U_{p}(B_{1})=U_{p}(B_{2})=1$. Hence
$B_{1}$ and $B_{2}$ are both abelian modulo their $O_{p'}$ subgroups.

Let $U=U_{\tilde q}(O_{p'}(B_{1}\cap B_{2}))$.
This group contains $U_{1}$ and is in particular nontrivial, and is $E_{1}$-invariant,
as well as
$U_{\tilde q}(O_{p'}(B_{1}))$ and $U_{\tilde q}(O_{p'}(B_{2}))$.
Now all these three subgroups are contained in a (unique) common
maximal $E_{1}$-invariant definable $\tilde q$-subgroup without elements
of order $p$ by the Uniqueness Theorem \ref{RelUniqLem}, say $\tilde U$.
Notice that $B_{1}=N\o(U_{\tilde q}(O_{p'}(B_{1})))$
and $B_{2}=N\o(U_{\tilde q}(O_{p'}(B_{2})))$ by maximality of
$B_{1}$ and $B_{2}$. Now applying the normalizer condition
\cite[Proposition 2.8]{FreconJaligot07} in the subgroup $\tilde U$ without elements of order $p$
yields easily
$U_{\tilde q}(O_{p'}(B_{1}))=\tilde U=U_{\tilde q}(O_{p'}(B_{2}))$.
Taking their common connected normalizers, $B_{1}=B_{2}$.

Our next claim follows from the same argument.

For the last claim, we note that there exists an element $s$ in $E_1$ of order $p$
such that $C_{U_{1}}(s)$ is nontrivial. By Lemma \ref{FaitDefGpAutpGp}
the latter is a $\tilde q$-group, and of course it is $E_{1}$-invariant. So
the preceding uniqueness applies to $C_{U_{1}}(s)$, and as
$C_{U_{1}}(s) \leq U_{1} \leq B$ we get that $B$ is the unique
maximal $E_{1}$-invariant definable connected solvable subgroup
containing $C_{U_{1}}(s)$.
But $C_{U_{1}}(s) \leq C\o(s)\leq B_{s}$ for some Borel subgroup
$B_{s}$ and $E_{1}\leq B_{s}$, so $B_{s}$ satisfies the same conditions as $B$,
so $B_{s}\leq B$ and $B=B_{s}$ is a Borel subgroup of $G$.

$(2)$. Suppose toward a contradiction
$U:=U_{\tilde q}(O_{p'}(C\o(E_{1})))$ nontrivial. It is of course
$E_{1}$-invariant. Recall that $Q$ is a fixed Carter
subgroup of $G$ containing the maximal $p$-torus $S$. As $Q\leq C\o(E_{1})$,
$Q$ normalizes the subgroup $U$. Now for any element $s$ of order $p$ in $E_{1}$
we have $UQ\leq C\o(s)$.

As $E_{1}\leq Q$, any Borel subgroup containing $UQ$ is $E_{1}$-invariant, and by
the first point there is a {\em unique} Borel subgroup containing $UQ$.
Now $C\o(s)$ is necessarily contained in this unique
Borel subgroup containing $UQ$, and this holds for any element $s$ of order
$p$ of $E_{1}$. We get a contradiction
to Lemma \ref{LemDispatchCentInvBorels}.
\qed

\bigskip
We note that the proof of the second point in Corollary \ref{CorUnipTheoAction}
actually shows that any definable connected subgroup containing $E_{1}$ and $U_{1}$
for some nontrivial $E_{1}$-invariant definable $\tilde q$-subgroup $U_{1}$
without elements of order $p$ is contained in a {\em unique} Borel subgroup
of $G$. Furthermore with the notation of Corollary \ref{CorUnipTheoAction} $(1)$
we have in any case
$N(U_{1})\cap N(E_{1})\leq N(U_{\tilde q}(O_{p'}(B)))=N(B)$.

\subsection{Bounding the Pr\"ufer rank via ``signalizer functors"}

There are two possible ways to prove that the Pr\"ufer $p$-rank is $2$ at this stage.
One may use the Uniqueness Theorem \ref{RelUniqLem} provided by
the $*$-local$_{\circ}\o$ solvability of the ambient group, or use the general
signalizer functor theory, which gives similar consequences in more general
contexts. We now explain how to use the signalizer functor theory to get
the bound on the Pr\"ufer $p$-rank, but we will rather continue the analysis
with the Uniqueness Theorem \ref{RelUniqLem} which is closer in spirit
to \cite[Lemma 6.1]{CherlinJaligot2004} and our original proof.
It also gives much more information in the specific context under consideration,
including when the Pr\"ufer $p$-rank is $2$, while the general
signalizer functor theory just provides the bound.

For $s$ a nontrivial element of $E$ we let
$$\theta(s)=U_{\tilde q}(O_{p'}(C(s))).$$
If $t$ is another nontrivial element
of $E$, then it normalizes the connected nilpotent $\tilde q$-group without
$p$-elements of order $\theta(s)$, and by
Lemmas \ref{FaitDefGpAutpGp} and \ref{FaitH/Op'HDivAbelian},
$C_{\theta(s)}(t)\leq U_{\tilde q}(O_{p'}(C(t)))=\theta(t)$. Hence
one has the two following properties:
\begin{itemize}
\item[$(1)$]
$\theta(s)^{g}=\theta(s^{g})$ for any $s$ in $E\setminus \{1\}$ and any $g$ in $G$.
\item[$(2)$]
$\theta(s) \cap C_{G}(t) \leq \theta(t)$ for any $s$ and $t$ in $E\setminus \{1\}$.
\end{itemize}
In the parlance of finite group theory one says that $\theta$ is an
{\em $E$-signalizer functor} on $G$. In groups of finite Morley rank
one says that $\theta$ is a {\em connected nilpotent}
$E$-signalizer functor, as any $\theta(s)$ is connected (by definition) and
nilpotent, which follows from Corollary \ref{CorEveryThingInFOp'}.
When $E_{1}$ is a subgroup of $E$ one defines
$$\theta(E_{1})=\<\theta(s)~|~s\in E_{1}\setminus \{1\}\>.$$

In groups of finite Morley rank there is no ``Solvable Signalizer Functor Theorem"
available as in the finite case \cite[Chapter 15]{MR1264416}
(see \cite{MR0323904, MR0297861, MR0417284, Bender75} for the history of the
finite case). However Borovik imported from finite group theory a
``Nilpotent Signalizer Functor Theorem" for groups of finite Morley rank
\cite{Borovik95} \cite[Theorem B.30]{BorovikNesin(Book)94},
stated as follows in \cite[Theorem A.2]{Burdges03} (and which suffices by the
unipotence theory of \cite{Burdges03} for which it has been designed
originally).

\begin{fact}
{\bf (Nilpotent Signalizer Functor Theorem)}
Let $G$ be a group of finite Morley rank, $p$ a prime, and $E\leq G$ a finite elementary
abelian $p$-group of order at least $p^{3}$. Let $\theta$ be a connected nilpotent
$E$-signalizer functor. Then $\theta(E)$ is nilpotent. Furthermore
$\theta(E)=O_{p'}(\theta(E))$ and $\theta(s)=C_{\theta(E)}(s)$ for any
$s$ in $E\setminus \{1\}$.
\end{fact}

(In the finite group theory terminology one says that $\theta$ is {\em complete}
when it satisfies the two properties of the last statement.)

In our situation one thus has, assuming toward a contradiction the Pr\"ufer $p$-rank
to be at least $3$, that $\theta(E)$ is nilpotent. Notice that the
definable connected subgroup $\theta(E)$ is nontrivial,
as $\theta(s)$ is nontrivial at least
for some $s$ by Lemma \ref{FaitDefGpAutpGp} and Fact \ref{BiGenerationWithpElts}.
In particular $N\o(\theta(E))$ is solvable by
$*$-local$_{\circ}\o$ solvability of $G$.

From this point on one can use arguments formally identical to those of
\cite[\S6.2-6.3]{Borovik95} used there for dealing with
``proper $2$-generated cores".

If $E_{1}$ and $E_{2}$ are two subgroups of $E$ of order at least $p^{2}$, then
for any $s$ in $E_{1}\setminus \{1\}$ one has
$\theta(s)\leq \<C_{\theta(s)}(t)~|~t\in E_{2}\setminus \{1\}\>\leq
\theta(E_{2})$ and thus $\theta(E_{1})=\theta(E_{2})$.

In particular $\theta(E)=\theta(E_{1})$ for any subgroup $E_{1}$ of $E$ of order at
least $p^{2}$.

Now if $g$ in $G$ normalizes such a subgroup $E_{1}$, then
$\theta(E)^{g}=\theta(E_{1})^{g}=\theta(E_{1}^{g})=
\theta(E_{1})=\theta(E)$ and thus $g\in N(\theta(E))$.

Take now as in Lemma \ref{LemDispatchCentInvBorels}
an element $s$ of order $p$ in $E$ such that $C\o(s)\nleq N\o(\theta(E))$.

Then, still assuming $E$ of order at least $p^{3}$, there exists a subgroup $E_{2}$ of $E$
of order at least $p^{2}$ and disjoint from $\<s\>$.
By Lemma \ref{Bipgeneration},
$$C\o(s)=\<C_{C\o(s)}(t)~|~t\in E_{2}\setminus \{1\}\>.$$
But now if $t$ is in $E_{2}$ as in the above equality, then $E_1:=\<s,t\>$ has order
$p^{2}$ as $E_{2}$ is disjoint from $\<s\>$, hence
$C_{C\o(s)}(t) \leq C(s,t)\leq N(\<s,t\>)=N(E_1)\leq N(\theta(E))$,
and this shows that
$C\o(s)\leq N\o(\theta(E))$. This is a contradiction, and as our only extra
assumption was that the Pr\"ufer $p$-rank was at least $3$, it must be $2$.

\subsection{Bounding the Pr\"ufer rank: Uniqueness Methods}

Anyway, we can get the bound similarly, by using more directly the
Uniqueness Theorem \ref{RelUniqLem} here instead of the axiomatized signalizer
functor machinery. Actually the proof below is the core of the proof of the
Nilpotent Signalizer Functor Theorem, and the Uniqueness Theorem here
shortcuts the passage to a quotient for the induction in the general case (see \cite{Burdges03}).

\begin{theo}
$S$ has Pr\"ufer $p$-rank $2$.
\end{theo}
\proof
Assume towards a contradiction $E$ has order at least $p^{3}$.

We then claim that there exists a {\em unique} maximal nontrivial $E$-invariant definable
$\tilde q$-subgroup without elements of order $p$.
Let $U_{1}$ and $U_{2}$ be two such subgroups. Then by
Lemma \ref{FaitDefGpAutpGp} and Fact \ref{BiGenerationWithpElts}
$C_{U_1}(E_1)$ and $C_{U_2}(E_2)$ are nontrivial $\tilde q$-subgroups for
some subgroups $E_1$ and $E_2$ of $E$, each of index $p$ in $E$. Assuming
$|E|\geq p^3$ then gives an element $s$ of order $p$ in $E_1 \cap E_2$.
Now $C_{U_1}(s)$ and $C_{U_2}(s)$ are nontrivial $\tilde q$-subgroups by Lemma \ref{FaitDefGpAutpGp}.
Clearly both are $E$-invariant, as $E$
centralizes $s$, and included in $U_{\tilde q}(O_{p'}(C\o(s)))$ as usual, which is also
$E$-invariant. Now the Uniqueness Theorem \ref{RelUniqLem} gives $U_{1}=U_{2}$,
as desired.

Hence there is a unique maximal $E$-invariant definable
$\tilde q$-subgroup without elements of order $p$, say ``$\theta(E)$" in the
notation of the signalizer functor theory. For the same reasons as mentioned
above, Lemma \ref{FaitDefGpAutpGp} and Fact \ref{BiGenerationWithpElts}, it is nontrivial.

Now by Lemma \ref{FaitDefGpAutpGp} and Fact \ref{BiGenerationWithpElts} again,
$C_{\theta(E)}(E_{1})$ is a nontrivial definable $\tilde q$-subgroup of
$\theta(E)$ for some subgroup $E_{1}$ of $E$ of index $p$.
As $U_{p}(C\o(E_{1}))=1$, the quotient $C\o(E_{1})/O_{p'}(C\o(E_{1}))$ is abelian
as usual, and the definable connected subgroup $C_{\theta(E)}(E_{1})$ is in
$O_{p'}(C\o(E_{1}))$, and in $U_{\tilde q}(O_{p'}(C\o(E_{1})))$.

But as $|E|\geq p^{3}$, $|E_{1}|\geq p^{2}$, and we get a contradiction to
Corollary \ref{CorUnipTheoAction} $(2)$.
\qed

\medskip
This proves clause $(2)$ of the Main Theorem %\ref{TheoPruferRankspStrongEmbed}
and completes its proof.
%of Theorem \ref{TheoPruferRankspStrongEmbed}
\qed

\section{Afterword}\label{S:Final}

We can also record informally some information gained along the proof of case
$(2)$ of our Main Theorem,
%\ref{TheoPruferRankspStrongEmbed}
which can be compared to \cite[6.1-6.6]{CherlinJaligot2004}. We let
$G$ and $S$ be as in case $(2)$ of the Main Theorem %\ref{TheoPruferRankspStrongEmbed},
and $Q$ be a Carter subgroup of $G$ containing $S$. Then $Q$ is contained in
at least two distinct Borel subgroups of $G$ by
Lemma \ref{LemDispatchCentInvBorels}, and in particular
$Q$ is divisible abelian by Fact \ref{corUnicityBorelSuperMax}
and \cite[Proposition 4.46]{DeloroJaligotI}. Now there are unipotence parameters
$\tilde q\neq (\infty,0)$ as in the proof of case $(2)$ of the Main Theorem
%\ref{TheoPruferRankspStrongEmbed}
(maybe one for $q=\infty$,
several for $q$ prime, but none for $q=p$ by Lemma \ref{LemAnyBorelContspPerp}).
All the results of the above analysis apply, now with $|\Omega_{1}(S)|=p^{2}$
necessarily.

By Corollary \ref{CorUnipTheoAction},
$$U_{\tilde q}(O_{p'}(C\o(\Omega_{1}(S))))=1.$$
As $\Omega_{1}(S)$ has order $p^{2}$, it contains in particular
$$\frac{p^{2}-1}{p-1}=p+1$$
pairwise noncollinear elements. It follows that there are at most
$p+1$ nontrivial subgroups of the form $U_{\tilde q}(O_{p'}(C\o(s)))$
for some nontrivial element $s$ of order $p$ of $S$, and at most $p+1$
Borel subgroups $B$ containing $Q$ (actually $\Omega_{1}(S)$-invariant suffices
as noticed after Corollary \ref{CorUnipTheoAction})
and such that $U_{\tilde q}(O_{p'}(B))\neq 1$.
By Corollary \ref{CorUnipTheoAction}, any such Borel subgroup
would contain $C\o(s)$ for any element
$s$ of order $p$ of $S$ having a nontrivial centralizer in $U_{\tilde q}(O_{p'}(B))$,
and $\Omega_{1}(S)$ has a trivial centralizer in
$U_{\tilde q}(O_{p'}(B))$.

The following corollary to the Main Theorem %\ref{TheoPruferRankspStrongEmbed}
will be of crucial use in \cite{DeloroJaligotII} to get a bound on Pr\"ufer $2$-ranks.

\begin{cor}\label{CorDichotomyHighPruferRk}
Let $G$ be a connected nonsolvable $*$-locally$_{\circ}\o$-solvable group of finite Morley
rank and of Pr\"ufer $p$-rank at least $2$ for some prime $p$, and fix a maximal $p$-torus
$S$ of $G$. Let $X$ be a maximal exceptional (finite) subgroup of $S$
(as in \cite[Lemma 3.29]{DeloroJaligotI}),
$\overline{H}=C\o(X)/X$, $\overline{K}$ a minimal definable
connected nonsolvable subgroup of $\overline{H}$ containing $\overline{S}$,
and let
$$\overline{B}=\<C\o_{\overline{K}}(\overline{s})~|~
\overline{s}\in \Omega_{1}(\overline{S}) \setminus \{ \overline{1} \}\>.$$
Then:
\begin{itemize}
\item[$(1)$]
either $\overline{B}<\overline{K}$, in which case $\overline{B}$ is a Borel subgroup
of $\overline{K}$; and if in addition $\overline{S}$ is a Sylow $p$-subgroup of $N_{N_{\overline{K}}(\overline{B})}(\overline{S})$, then $N_{\overline{K}}(\overline{B})$
is $p$-strongly embedded in $\overline{K}$,
\item[$(2)$]
or $\overline{B}=\overline{K}$, in which case
$\overline{S}$, as well as $S$, has Pr\"{u}fer $p$-rank $2$.
\end{itemize}
\end{cor}
\proof
It suffices to apply our Main Theorem %\ref{TheoPruferRankspStrongEmbed}
in $\overline{K}$. We note that $\overline{S}$ and $S$ have the same
Pr\"ufer $p$-rank, as $X$ is finite by \cite[Lemma 3.18]{ DeloroJaligotI}.
\qed

\bigskip
Cases $(1)$ and $(2)$ of the Main Theorem %\ref{TheoPruferRankspStrongEmbed}
and
Corollary \ref{CorDichotomyHighPruferRk} correspond respectively to
Sections 7 and 6 of \cite{CherlinJaligot2004} in presence of divisible torsion.
%The remaining analysis of both of these sections, as well as the treatment without the extra assumption for $p$-strong embedding in case $(1)$, will be considered in our separate paper on Weyl groups, mentioned already in Section \ref{SectionLocAnalandGen}.

For $p=2$ case $(1)$ will entirely disappear in \cite{DeloroJaligotII} by an
argument similar to the one used in \cite[Case I]{BurdgesCherlinJaligot07} for minimal connected simple groups.

\bibliographystyle{alpha}
\bibliography{biblio}

\end{document}